\documentclass[english, 10pt]{article}

\usepackage{std_english}
\usepackage{subfig}

\numberwithin{equation}{section}

\newcommand{\wght}[1]{\zset_+^{\vts{#1}}}
\newcommand{\ssp}[1]{\mathrm{STAB}(#1)}

\newcommand{\ecol}{\Lambda}

\newcommand{\fcn}[2]{\chi_f(#2,#1)}
\newcommand{\cn}[2]{\chi(#2,#1)}
\newcommand{\sts}[1]{\mathcal{S}(#1)}

\newcommand{\cfpfree}{odd-$C_5^{+}$-free\xspace}

\newcommand{\qset}{\mathbb{Q}}

\newcommand{\psbg}[2]{#1\cro{#2}}

\declaretheorem[title=Theorem]{thm}

\declaretheorem[title=Corollary]{cor}

\declaretheorem[title=Proposition]{prop}

\declaretheorem[title=Lemma]{lem}

\theorembodyfont{\rmfamily}

\declaretheorem[title=Proof,
				postheadhook=--~, 
				numbered=no,
				postfoothook=\hfill$\blacksquare$]{dem}
				
\theoremheaderfont{\mdseries\itshape}

\declaretheorem[title=Case,
				numberwithin=dem,
				]{case}

\newcommand{\vts}[1]{V(#1)}
\newcommand{\eds}[1]{E(#1)}

\newcommand{\vn}{\varnothing}

\newcommand{\se}{\subseteq}

\newcommand{\sm}{\setminus}

\newcommand{\cro}[1]{\left[#1\right]}

\newcommand{\ceil}[1]{\lceil #1\rceil}

\newcommand{\rset}{\mathbb{R}}
\newcommand{\zset}{\mathbb{Z}}

\newcommand{\set}[1]{\left\{#1\right\}}

\newcommand{\myrefwp}[1]{\cref{#1}\xspace}
\title{Integer round-up property for the chromatic number of some h-perfect graphs}
\author{Yohann \sc{Benchetrit}
\thanks{CNRS, G-SCOP, F-38000 Grenoble, France \newline
Univ. Grenoble Alpes, G-SCOP, F-38000 Grenoble, France}}

\begin{document}
\maketitle

\begin{abstract}
A graph is h-perfect if its stable set polytope can be completely described by non-negativity, clique and odd-hole constraints. It is  t-perfect if it furthermore has no clique of size 4.

For every graph $G$ and every $c\in\wght{G}$, the weighted chromatic number of $(G,c)$ is the minimum cardinality of a multi-set $\mathcal{F}$ of stable sets of $G$ such that every $v\in\vts{G}$ belongs to at least $c_v$ members of $\mathcal{F}$.

We prove that every h-perfect line-graph and every t-perfect claw-free graph $G$ has the integer round-up property for the chromatic number: for every non-negative integer weight $c$ on the vertices of $G$, the weighted chromatic number of $(G,c)$ can be obtained by rounding up its fractional relaxation. In other words, the stable set polytope of $G$ has the integer decomposition property. 

Another occurrence of this property was recently obtained by Eisenbrand and Niemeier for fuzzy circular interval graphs (extending previous results of Niessen, Kind and Gijswijt). These graphs form another proper subclass of claw-free graphs. 

Our results imply the existence of a polynomial-time algorithm which computes the weighted chromatic number of t-perfect claw-free graphs and h-perfect line-graphs. Finally, they yield a new case of a conjecture of Goldberg and Seymour on edge-colorings.

\end{abstract}

\section{Introduction}
In this paper, we consider undirected finite graphs. They can have multiple edges but no loops. We say that a graph is \emph{simple} if it does not have two parallel edges.

Let $G$ be a graph. We write $\vts{G}$ for the vertex set of $G$ and $\eds{G}$ for its edge set.
A set $S\se\vts{G}$ is \emph{stable} if its elements are pairwise non-adjacent in $G$. The set of stable sets of $G$ is denoted by $\sts{G}$.

For every $c\in\wght{G}$, the \emph{weighted chromatic number of $(G,c)$}, denoted $\cn{c}{G}$, is the minimum cardinality of a multi-set $\mathcal{F}$ (counting multiple occurrences of the same element) of stable sets of $G$ such that every $v\in\vts{G}$ belongs to at least $c_v$ members of $\mathcal{F}$.
In other words:
\begin{equation}\label{eqn-ip-form}
 \cn{c}{G}=\min\set{
\sum_{S\in\sts{G}}y_S\colon\,\, 
y\in\zset_+^{\sts{G}};\,
\sum_{S\in\sts{G}\colon v\in S} y_S\geq c_v,\,\forall v\in\vts{G}
}.
\end{equation}
The \emph{chromatic number} $\chi(G)$ is equal to $\cn{\mathbf{1}}{G}$, where $\mathbf{1}$ is the all-1 vector of $\zset^{\vts{G}}$. We will speak of the \emph{unweighted case} when considering the weight function $\mathbf{1}$. 

Replacing $\zset$ with $\qset$ in \eqref{eqn-ip-form}, we obtain a linear program whose optimum value is the \emph{weighted fractional chromatic number of $(G,c)$}. We write it $\fcn{c}{G}$ (and simply $\chi_f(G)$ in the unweighted case). Hence, the inequality $\ceil{\fcn{c}{G}}\leq\cn{c}{G}$ always holds.

The chromatic number of a graph has been extensively studied in various contexts of discrete optimization and graph theory (see for example the book \cite{Jensen1995}). Karp \cite{Karp1972} proved that it is NP-hard to compute and several inapproximability results  were later obtained (see Huang \cite{Huang2013} for a recent example). 
Finding its fractional counterpart is also an NP-hard problem in general since it is equivalent to the maximum-weight stable set problem, through the ellipsoid method \cite{Groetschel1988} (a proof of this fact using a Karp-reduction is not known).

A graph $G$ is \emph{perfect} if every induced subgraph $H$ of $G$ satisfies $\chi(H)=\omega(H)$, where $\omega(H)$ denotes the maximum cardinality of a clique of $H$.
Results of Lov\'asz \cite{Lovasz1972} and Fulkerson \cite{Fulkerson1972} imply, as stated by Chv\'atal \cite{Chvatal1975}:
\begin{thm}[Lov\'asz, Fulkerson, Chv\'atal]\label{thm-cfl-rep}
For every perfect graph $G$ and every $c\in\zset_+^{\vts{G}}$: \[ \cn{c}{G}=\fcn{c}{G}.  \]
\end{thm}
It follows from results of Gr\"{o}tschel, Lov\'asz, Schrijver \cite{Groetschel1988} that \emph{the (integer or fractional) weighted chromatic number of a perfect graph can be found in polynomial-time}. 

A graph is \emph{h-perfect} if the non-trivial facets of its stable set polytope are given by non-negativity, clique and odd-hole constraints. 
H-perfect graphs form a (proper) superclass of perfect graphs. By Gr\"{o}tschel, Lov\'asz, Schrijver \cite{Grotschel1986}, \emph{the weighted fractional chromatic number of an h-perfect graph can also be computed in polynomial-time}. 

However, the complexity of determining $\cn{c}{G}$ in h-perfect graphs is unknown. Hence, the study of the gap between the weighted chromatic number of $(G,c)$ and its fractional version may help to design (either exact or approximation) polynomial-time algorithms for this problem. We do not know of any result giving a bound on this gap for every h-perfect graph and every weight.

A graph $G$ has the \emph{integer round-up property} (for the chromatic number) if $\cn{c}{G}=\ceil{\fcn{c}{G}}$ holds for every $c\in\zset_+^{\vts{G}}$.
Therefore, if every graph of a subclass $\mathcal{G}$ of h-perfect graphs has this property, then their weighted chromatic number can be computed  in polynomial-time (for every weight).

The \emph{line graph} of a graph $G$, denoted $L(G)$, is the simple graph whose vertex set is $\eds{G}$ and whose edge set is formed by the pairs of incident edges of $G$. A graph is \emph{claw-free} if it does not contain the complete bipartite graph $K_{1,3}$ as an induced subgraph. Claw-free graphs form a proper superclass of line graphs.

The class of h-perfect claw-free graphs was investigated by Bruhn and Stein in \cite{Bruhn2012}. In particular, they proved the unweighted case of the integer round-up property for these graphs:

\begin{thm}[Bruhn, Stein \cite{Bruhn2012}]\label{thm-brst-ir}
Every h-perfect claw-free graph $G$ satisfies $\chi(G)=\ceil{\chi_f(G)}$.
\end{thm}

The line graph of the Petersen graph shows that this is not true for line graphs in general.
In this paper, we extend this result to arbitrary weights for two subclasses of h-perfect claw-free graphs. 
First, we will show that:

\begin{thm}\label{thm-ben-line-irp}
Every h-perfect line-graph has the integer round-up property.
\end{thm}

The corresponding unweighted result was obtained by Bruhn and Stein \cite{Bruhn2012} and serves as a lemma for \myrefwp{thm-brst-ir}. The proof consists in coloring the edges of the source graph, whose structure is described by a theorem of Cao and Nemhauser \cite{Cao1998}.We follow the same idea to show \myrefwp{thm-ben-line-irp}. Further arguments are needed to handle phenomena which occur only in the weighted case.

The \emph{chromatic index} of a graph $G$, denoted $\chi'(G)$, is the minimum number of colors needed to assign to each edge of $G$ a color such that two incident edges receive different colors. 
Conjectures of Goldberg \cite{Goldberg1973} and Seymour \cite{Seymour1979a} imply that: \emph{for every graph $G$, $\chi'(G)$ is equal to 
$\ceil{\chi_f(L(G))}$ 
or $ \ceil{\chi_f(L(G))}+1 $}.
There are not many known classes of graphs that are defined by an excluded-subgraph (or minor) assumption and for which $\chi'(G)=\ceil{\chi_f(L(G))}$ holds for every member of the class. Seymour \cite{Seymour1990} showed that this holds for graphs without a subdivision of the complete graph $K_4$ (that is series-parallel graphs; Fernandes and Thomas \cite{Fernandes2013} later found a shorter proof) and Marcotte \cite{Marcotte1986} proved it for graphs which do not have a minor isomorphic to $K_5$ minus an edge. 
In Section 3, we prove (as an intermediate result towards theorem \ref{thm-ben-line-irp}) that every graph $G$ which does not contain (as a subgraph) a totally odd subdivision of $C_5^{+}$ (the graph obtained from the circuit of length 5 by adding an edge between two non-adjacent vertices) satisfies:  $\chi'(G)=\ceil{\chi_f'(G)}$ (see theorem \ref{thm-ecol-form}). 
These results do not imply one another.

The other main result of this paper concerns t-perfect graphs.
A graph is \emph{t-perfect} if it is h-perfect and does not contain a clique of size 4. 

\begin{thm}\label{thm-ben-cft-irp}
Every t-perfect claw-free graph has the integer round-up property.
\end{thm}

The unweighted case is a special case of theorem \ref{thm-brst-ir} and is obtained by a reduction to the line-graph case. We had to follow a new approach in proving \myrefwp{thm-ben-cft-irp}: if $G$ is a t-perfect claw-free  graph and $c\in\wght{G}$, then we can either reduce the size of $c$ using certain subgraphs (and use induction) or apply theorem \ref{thm-ben-line-irp}. 

Our proof strongly relies on the characterization of t-perfect squares of circuits, which was stated and proved in \cite{Bruhn2012}. Our starting point is Harary's characterization of line graphs (\cite{Harary1972} pg. 74: theorem 8.4), similarly as in the proof of the unweighted case. Then, we need merely the non-t-perfection of certain graphs \cite{Bruhn2012}.

To our knowledge, there are only two previous results on the integer round-up property for h-perfect (non-perfect) graphs: Kilakos and Marcotte \cite{Kilakos1997} proved it for graphs without a subdivision of $K_4$ and Gerards (unpublished, \cite{Schrijver2003} pg. 1207) later extended this result to graphs not containing an odd subdivision of $K_4$.
In \cite{Kilakos1997}, Kilakos and Marcotte developed a general method to prove the integer round-up property. We do not see how to apply this method to t-perfect claw-free graphs.

We do not know if every h-perfect claw-free graph has the integer round-up property. 
It does not hold for h-perfect graphs in general. Indeed, 
Shepherd \cite{Shepherd1995} noted that \emph{every t-perfect graph which has the property must have chromatic number at most 3}. However, Laurent and Seymour (\cite{Schrijver2003}, pg. 1207) found a t-perfect graph which is not 3-colorable.
Whether every 3-colorable t-perfect graph has the integer round-up property is still an open problem (which was already mentioned by Shepherd in \cite{Shepherd1995}). Theorem \ref{thm-ben-cft-irp} shows that this is at least the case for claw-free graphs.

Computing the chromatic number of an h-perfect graph (in the unweighted case) remains unsolved in general and no upper bound has been found.
Seb\H{o} conjectures that \emph{triangle-free t-perfect graphs are 3-colorable} and showed that, if valid, this would imply that every h-perfect graph $G$ is $(\omega(G)+1)$-colorable (see \cite{Bruhn2012}).

T-perfect graphs were introduced by Chv\'atal \cite{Chvatal1975} and h-perfect graphs were later defined by Sbihi and Uhry \cite{Sbihi1984}.
Examples of h-perfect graphs include almost-bipartite graphs \cite{Fonlupt1982} and more general classes of graphs defined by excluding (as a subgraph) non-t-perfect subdivisions of $K_4$ (see \cite{Boulala1979,Gerards1986,Gerards1998}). Shepherd \cite{Shepherd1995} characterized t-perfection among complements of line graphs. Bruhn and Stein \cite{Bruhn2012} recently gave a characterization of t-perfect claw-free graphs in terms of certain minor-operations (namely, t-minors) defined by Gerards and Shepherd \cite{Gerards1998}.

The maximum-weight stable set problem can be solved in polynomial-time in the class of h-perfect graphs \cite{Grotschel1986}. Eisenbrand et al. \cite{Eisenbrand2003} gave a combinatorial algorithm for the cardinality case in t-perfect graphs.
By contrast, the computational complexity of deciding t-perfection is unknown.

The integer round-up property (for the chromatic number) of $G$ corresponds to the \emph{integer decomposition property} of the stable set polytope $\ssp{G}$: for every $k\geq 1$, every integer vector of $k\cdot\ssp{G}$ is the sum of $k$ integer vectors of $\ssp{G}$. The integer decomposition property of polyhedra was introduced by Baum and Trotter \cite{Baum1981} as a general framework to study linear programs for which rounding (up or down) the optimum value yields the optimum value of the associated integer program.
It occurs in various contexts such as the cutting stock problem \cite{Scheithauer1997} and commutative algebra (under the notion of normality, see for example \cite{Ohsugi1998}).

Circular arc graphs form another class of claw-free graphs which have the integer round-up property. This was obtained by Niessen \cite{Niessen2000}, Gijswijt \cite{Gijswijt2005} and later extended to fuzzy circular interval graphs by Eisenbrand et al \cite{Eisenbrand2012} (both these classes are incomparable with the class of h-perfect claw-free graph in terms of inclusion). These graphs appear in the context of the problem of finding a nice description of the stable set polytope of claw-free graphs.

King and Reed \cite{King2013} proved that $\chi$ and $\chi_f$ agree asymptotically in quasi-line graphs (they form a proper subclass of claw-free graphs and superclass of line graphs). 
It is not difficult to show that every h-perfect claw-free graph is quasi-line, but we do not use it in our proofs.

We end this introduction with an outline of the paper.
In Section 2, we give definitions and basic properties for h-perfect graphs, weighted colorings and related notions.
Section 3 contains the proof of theorem \ref{thm-ben-line-irp}. We use it to prove theorem \ref{thm-ben-cft-irp} in Section 4.
In Section 5, we derive an explicit formula for the weighted chromatic number of h-perfect line-graphs and t-perfect claw-free graphs as a consequence of theorems \ref{thm-ben-line-irp} and \ref{thm-ben-cft-irp}. We also state a related formula for the chromatic index of \cfpfree graphs (see definition in Section 3).
Finally, the algorithmic aspects of our results are discussed in Section 6.

\section{Basic definitions and properties}

Let $G$ be a graph. 
A \emph{subgraph} of $G$ is a graph $H$ such that $\vts{H}\se\vts{G}$ and 
$\eds{H}\se\eds{G}$. A subgraph $H$ of $G$ is \emph{induced} if $\eds{H}$ contains every edge of $G$ whose ends both belong to $\vts{H}$. For every $U\se\vts{G}$, we write $G\cro{U}$ for the induced subgraph of $G$ whose vertex set is $U$.
The \emph{components} of $G$ are its maximal connected subgraphs (with respect to the subgraph relation). We say that $G$ is \emph{$k$-regular} (for a non-negative integer $k$) if each vertex of $G$ has exactly $k$ neighbors in $G$.

A \emph{clique} of $G$ is a set of pairwise-adjacent vertices of $G$. The \emph{triangles} of $G$ are its cliques of size 3.
A \emph{circuit} of $G$ is a 2-regular connected subgraph of $G$. It is \emph{odd} if it has an odd number of vertices. 
An \emph{odd hole} of $G$ is an induced odd-circuit of $G$ with at least 5 vertices. 

Let $X\se\vts{G}$. The \emph{characteristic vector} of $X$, denoted $\chi^{X}$, is the $\set{0,1}$-vector of $\qset^{\vts{G}}$ defined by: for every $v\in X$, $\chi^{X}(v)=1$ if and only if $v\in X$. If $X$ has only one element $v$, then its characteristic vector is written $\chi^{v}$.
The \emph{stable set polytope} of $G$, denoted $\ssp{G}$, is the convex hull of the characteristic vectors of the stable sets of $G$.

Since it is a full-dimensional polyhedron, $\ssp{G}$ can be uniquely described (up to non-zero scalar multiplication of inequalities) as the set of solutions of a system of linear inequalities. 
It is not possible to test in polynomial-time if a given vector belongs to $\ssp{G}$, unless $P=NP$.

Classes of valid inequalities for $\ssp{G}$ include the \emph{non-negativity inequalities} $x_v\geq 0$ (for every $v\in\vts{G})$ and the \emph{clique inequalities} $\sum_{v\in K}x_v\leq 1$ for every clique $K$ of $G$. 
Chv\'atal \cite{Chvatal1975} proved that \emph{these inequalities are enough to describe $\ssp{G}$ if and only if $G$ is perfect}.

Another important class of valid inequalities for $\ssp{G}$ is given by \emph{odd-hole inequalities} $\sum_{v\in\vts{C}}x_v\leq\frac{|\vts{C}|-1}{2}$, where $C$ is an odd hole of $G$.
A graph $G$ is \emph{h-perfect} if the non-negativity, clique and odd-hole inequalities are enough to describe $\ssp{G}$, in other words if:

\begin{equation}\label{eqn-hparf-def}
\ssp{G}=\left\{x\in\rset^{\vts{G}}\colon\begin{array}{cc}
  x\geq 0,  &  \\ 
  \displaystyle\sum_{v\in K}x_v\leq 1  & \text{$\forall K$ clique of $G$,}  \\ 
  \displaystyle\sum_{v\in \vts{C}}x_v\leq \dfrac{|\vts{C}|-1}{2}  & \text{$\forall C$ odd hole of $G$.} \\ 
\end{array}\right\}.
\end{equation}
Notice that this is equivalent to asking that the polyhedron described by the right side of this equality is \emph{integral} (that is, each of its vertices is an integer vector).
Since no induced subgraph of a perfect graph is an odd hole, perfect graphs are h-perfect. On the other hand, odd holes form the most basic example of non-perfect h-perfect graphs.
A graph is \emph{t-perfect} if it is h-perfect and does not contain a clique of size 4. 

\begin{figure}[h]
\centering
\includegraphics{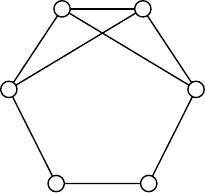}
\caption{Duplicating a vertex in the circuit of length 5}\label{fig-rep}
\end{figure}

Let $c\in\wght{G}$. A \emph{$k$-coloring} of $G$ (where $k$ is a non-negative integer) is a multi-set $\mathcal{F}$ of $k$ stable sets of $G$ such that every $v\in\vts{G}$ belongs to at least $c_v$ members of $\mathcal{F}$ (equivalently, to exactly $c_v$ members of $\mathcal{F}$). We say that $(G,c)$ is \emph{$k$-colorable} if it admits a $k$-coloring. 
Hence, $\cn{c}{G}$ is equal to the smallest integer $k$ such that $(G,c)$ is $k$-colorable.
Recall that $\sts{G}$ denotes the set of stable sets of $G$.
The \emph{weighted fractional chromatic number of $(G,c)$}, denoted $\chi_f(G,c)$, is the optimum value of the linear relaxation of the integer programming formulation \eqref{eqn-ip-form} for $\chi(G,c)$. That is:

\begin{equation}\label{eqn-lp-relax}
 \fcn{c}{G}=\min\set{
\sum_{S\in\sts{G}}y_S\colon\,\, 
y\in\qset_+^{\sts{G}};\,
\sum_{S\in\sts{G}\colon v\in S} y_S\geq c_v,\,\forall v\in\vts{G}
}. 
\end{equation}
When $c$ is the all-1 vector, we simply write it $\chi_f(G)$. 
Weighted and unweighted versions are related through the notion of \emph{replication}. For every $c\in\wght{G}$, let $G^{c}$ denote the graph obtained as follows: replace each vertex $v$ of $G$ with a clique $K_v$ of size $c_v$ (in particular, delete $v$ if $c_v=0$), and for every $u,v\in\vts{G}$: put every edge between two cliques $K_u,K_v$ if $uv\in\eds{G}$ and no edge otherwise. It is straightforward to check the following relations:

\begin{prop}\label{prop-chi-repli}
For every graph $G$ and every $c\in\wght{G}$: 
\begin{center}
$\chi(G,c)=\chi(G^{c})$ and $\chi_f(G,c)=\chi_f(G^{c})$.
\end{center}
\end{prop}
An important step in the proof of theorem \ref{thm-cfl-rep} is the \emph{replication lemma} of Lov\'asz \cite{Lovasz1972}: \emph{for every perfect graph $G$ and every $c\in\wght{G}$, the graph $G^{c}$ is perfect. }
The replication lemma is not true for t-perfect claw-free graphs. Indeed, the circuit of length 5 is t-perfect and duplicating one of its vertices (see figure \ref{fig-rep}) gives a non-t-perfect graph (see lemma 6 for a simple argument). 

In Section 4, we will use the following basic property of h-perfection (which is obtained by an easy polyhedral argument):

\begin{prop}\label{prop-tperf-indu}
Every induced subgraph of an h-perfect graph is h-perfect. In particular, t-perfection is closed under taking induced subgraphs.
\end{prop}
For every graph $G$ and $c\in\wght{G}$, let
$  \omega(G,c)=\max\set{   \sum_{v\in K}c_v   \colon  \text{$K$ clique of $G$}    } $ and:
\[  \Gamma(G,c)=\max\set{   \frac{2}{|\vts{C}|-1}\sum_{v\in\vts{C}}c_v \colon \text{$C$ odd hole of $G$}    }.  \] 
Using the duality theorem of linear programming on \ref{eqn-lp-relax}, it is straightforward to show the following formula:

\begin{prop}\label{prop-fcn-hparf}
For every h-perfect graph $G$ and every $c\in\wght{G}$:
\[ \chi_f(G,c)=\max(\omega(G,c),\Gamma(G,c)). \]
\end{prop}

\section{H-perfect line-graphs}
The purpose of this section is to prove theorem \ref{thm-ben-line-irp}.
In part 3.1, we state an edge-coloring result (theorem \ref{thm-ecol-form}) and show that it easily implies theorem \ref{thm-ben-line-irp} using a result of Cao and Nemhauser (which is a direct consequence of Edmonds' description of the matching polytope).

Part 3.2 is devoted to the proof of this edge-coloring statement. It relies on an auxiliary result (lemma \ref{thm-ben-crit}) whose proof is postponed to part 3.3.

\subsection{Reduction to an edge-coloring result} 
Let $H$ be a graph. 
Recall that the \emph{line graph of $H$} is denoted $L(H)$ and is the simple graph whose vertex set is $\eds{H}$ and whose edge set is formed by the pairs of incident edges of $H$.
We say that $H$ is a \emph{line graph} if there exists a graph $H'$ such that $H=L(H')$. 

A \emph{path} is a graph obtained from a circuit by deleting an edge. The \emph{length} of a path $P$ is the cardinality of its edge set. We say that $P$ is \emph{odd} (resp. \emph{even}) if it has odd (resp. \emph{even}) length.

A graph $K$ is a \emph{subdivision of $H$} if it is obtained from $H$ by replacing every edge $e\in\eds{H}$ with a path $P_e$ (of non-zero length) joining the ends of $e$ such that: for every pair of edges $e,f\in\eds{H}$, the paths $P_e$ and $P_f$ do not share inner vertices. If each $P_e$ has odd length, then we furthermore say that $K$ is a \emph{totally odd subdivision of $H$}.

\begin{figure}[h]
\centering
\includegraphics{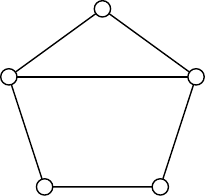}
\caption{The graph $C_5^{+}$, also known as the \emph{house}}\label{fig-c5p}
\end{figure}

Let $C_5^{+}$ be the graph defined by $\vts{C_5^{+}}=\set{1,2,3,4,5}$ and $\eds{C_5^{+}}=\set{12,23,34,45,15,25}$ (see figure \ref{fig-c5p}). An \emph{odd-$C_5^{+}$ of $H$} is a subgraph of $H$ which is isomorphic to a totally odd subdivision of $C_5^{+}$. A graph is \emph{odd-$C_5^{+}$-free} if it does not contain an odd-$C_5^{+}$ as a subgraph.

Various notions of perfection in line graphs were investigated by Cao and Nemhauser. In particular, they characterized the h-perfection of $L(H)$ in terms of the exclusion of certain graphs as subgraphs of $H$.
We need only part of their result:

\begin{thm}[Cao, Nemhauser \cite{Cao1998}]\label{thm-cao-eq}
Let $H$ be a graph. 
If $L(H)$ is h-perfect, then $H$ is \cfpfree.
\end{thm}
Let us recall the basic terminology of edge-coloring.
A \emph{matching} of a graph $H$ is a set of pairwise non-incident edges of $H$. A \emph{$k$-edge-coloring} (where $k$ is a non-negative integer) is a set $\ecol$ of $k$ matchings of $H$ such that every edge belongs to at least one element of $\ecol$ (equivalently, to exactly one element of $\ecol$). 
The \emph{chromatic index of $H$}, denoted $\chi'(H)$, is the least integer $k$ such that $G$ admits a $k$-edge-coloring. In other words: $\chi'(H)=\chi(L(H))$. The \emph{fractional chromatic index of $H$} is analogously defined by $\chi_f'(H)=\chi_f(L(H))$. 

\begin{thm}\label{thm-ecol-form}
Every \cfpfree graph $H$ satisfies: $\chi'(H)=\ceil{\chi_f'(H)}$.
\end{thm}

We first show that this result and theorem \ref{thm-cao-eq} imply theorem \ref{thm-ben-line-irp}.

\begin{dem}[of theorem \ref{thm-ben-line-irp}]
Let $G$ be an h-perfect line-graph, $c\in\wght{G}$ and $H$ be a graph such that $G=L(H)$. Let $H'$ be the graph obtained from $H$ by replacing each edge $e=uv\in\eds{H}$ by $c_e$ parallel edges between $u$ and $v$. It is straightforward to check that $G^{c}=L(H')$.
By proposition \ref{prop-chi-repli}, $\cn{c}{G}=\chi(G^{c})=\chi'(H')$ and similarly, $\fcn{c}{G}=\chi_f'(H')$. 

By theorem \ref{thm-cao-eq}, $H$ is \cfpfree. Therefore, $H'$ is also \cfpfree (it only depends on its underlying simple graph, which is the same as $H$) and the conclusion follows from theorem \ref{thm-ecol-form}.

\end{dem}

\subsection{Proof of theorem \ref{thm-ecol-form}}
Let $H$ be a graph. An \emph{odd ring} of a graph $H$ is an induced subgraph of $H$ whose underlying simple graph is an odd circuit (which can be of length 3).
Let:
\[ \Gamma '(H)=\max\set{ \frac{2}{|\vts{R}|-1}|\eds{R}|\colon \text{$R$ is an odd ring of $H$}              }, \]
The \emph{degree} of a vertex $v$ of $H$, denoted $d_H(v)$, is the number of edges of $H$ incident with $v$ (since $H$ can have multiple edges, $d_H(v)$ can be different from the number of neighbors of $v$ in $H$). 
Let $\Delta(H)$ denote the largest degree of a vertex of $H$.

Let $e\in\eds{H}$. The graph $H-e$ is defined by $\vts{H-e}=\vts{H}$ and $\eds{H-e}=\eds{H}\sm\set{e}$ (notice that the other edges parallel to $e$ (if any) are not deleted). 
The edge $e$ is \emph{critical} if $\chi'(H-e)<\chi'(H)$ (that is, $\chi'(H-e)=\chi'(H)-1$). 

The main ingredient of the proof of theorem \ref{thm-ecol-form} is the following "concentration" lemma (whose proof is postponed to part 3.3): 
\begin{lem}\label{thm-ben-crit}
Let $H$ be a graph such that $\chi'(H)>\Delta(H)$ and let $e\in\eds{H}$.
If $e$ is critical and is not an edge of an odd-$C_5^{+}$ of $H$, then there exists an odd ring $R$ of $H$ such that $e\in\eds{R}$ and:  
\[ |\eds{R}|= r\cdot\chi'(H-e)+1, \]
where $r=\frac{|\vts{R}|-1}{2}$. 

\end{lem}

We need one more result on the fractional chromatic index of a graph. In section 5, we will see that equality actually holds in the statement i) below for every \cfpfree graph. As a byproduct, we will obtain a new formula for the chromatic index of these graphs. We defer this formula for the sake of clarity, since only the lower bound is needed in the proof of theorem \ref{thm-ecol-form}.

\begin{prop}\label{prop-basic-propty}
Let $H$ be a graph. The following statements hold:
\begin{itemize}
\item[i)] $\chi_f'(H)\geq \max(\Delta(H),\Gamma'(H))$,
\item[ii)] for every subgraph $K$ of $H$, we have $\chi_f'(K)\leq\chi_f'(H)$.
\end{itemize}
\end{prop}
\begin{dem}
Let $\mathcal{M}(H)$ be the set of all matchings of $H$. By the duality theorem of linear programming, we have:
\[ \chi_f'(G)=\max\set{\sum_{e\in\eds{H}}x_e
\colon x\in\qset_+^{\eds{G}};
\, \sum_{e\in M}x_e\leq 1,\,\text{for every $M\in\mathcal{M}(H)$}}. \]
Let $v\in\vts{H}$ and $R$ be an odd ring of $H$. Let $\delta(v)$ denote the set of edges incident with $v$.
Clearly, each matching of $H$ contains at most one edge of $\delta(v)$ and at most $\frac{2}{|\vts{R}|-1}|\eds{R}|$ edges of $R$. Hence, both $\chi^{\delta(v)}$ and $\frac{2}{|\vts{R}|-1}\chi^{\eds{R}}$ are feasible solutions of the linear program above and this implies i).

Statement ii) follows from the fact that any optimal solution $x$ of the linear program above for $\chi_f'(K)$ can be extended to a feasible solution of the program for $\chi_f'(H)$, by setting $x_e=0$ for every $e\in\eds{H}\setminus\eds{K}$. 

\end{dem}

We are now ready to prove theorem \ref{thm-ecol-form} (the proof of lemma \ref{thm-ben-crit} being the subject of the next part):

\begin{dem}[of theorem \ref{thm-ecol-form}]
For every graph $G$, let $\kappa(G)$ denote $\ceil{\chi_f'(G)}$.

Seeking a contradiction, let $H$ be an \cfpfree graph with $\chi'(H)\neq\kappa(H)$ and choose $|\eds{H}|$ minimum. We actually have $\chi'(H)>\kappa(H)$ (since $\chi'(G)\geq\kappa(G)$ clearly holds for every graph $G$).

By proposition \ref{prop-basic-propty}, $\chi'(H)>\Delta(H)$.
Let $e\in\eds{H}$ and put $H'=H-e$. We have $\chi'(H')=\kappa(H')$ because of the minimality of $H$.
Besides,  proposition \ref{prop-basic-propty}.ii) implies that $\kappa(H')\leq\kappa(H)$, so $e$ must be a critical edge of $H$.
Since $H$ is \cfpfree, lemma \ref{thm-ben-crit} can be applied to $H$ and $e$. 

Hence, $H$ has an odd ring $R$ such that $e\in\eds{R}$ and $|\eds{R}|=r\cdot\chi'(H')+1$, where $r=\frac{|\vts{R}|-1}{2}$. 
By proposition \ref{prop-basic-propty}.i), $\kappa(H)\geq \frac{|\eds{R}|}{r}$. 
Therefore, 
$\kappa(H)>\chi'(H')$ and this contradicts our assumption on $H$.
\end{dem}

\subsection{Proof of lemma \ref{thm-ben-crit}}
Let $H$ be a graph. The \emph{multiplicity} of an edge $e$ of $H$, denoted $\mu_H(e)$, is the number of edges which are parallel to $e$ (including $e$). 

Let $v$ be a vertex and $M$ be a matching of $H$. We say that $M$ \emph{covers} $v$ if $M$ has an edge incident with $v$, and that it \emph{misses} $v$ otherwise. An edge-coloring of $H$ is \emph{optimal} if it uses exactly $\chi'(H)$ colors. 

\begin{prop}\label{lem-precol}
Let $H$ be a graph such that $\chi'(H)>\Delta(H)$ and let $e=uv$ be a critical edge of $H$. If $\ecol$ is an optimal edge-coloring of $H-e$, then:
\begin{itemize}
\item [i)] every matching $M\in\ecol$ covers at least one of $u$ and $v$ ,
\item [ii)] there exist two matchings $A,B\in \ecol$ such that $A$ covers $u$ and misses $v$, whereas $B$ covers $v$ and misses $u$.
\end{itemize}
\end{prop}
\begin{dem}
If i) did not hold, then  $\ecol$ could be extended to an edge-coloring $\ecol'$ of $H$ (by adding $e$ to a matching which misses both $u$ and $v$). This would contradict $\chi'(H-e)<\chi(H)$.

We now prove ii). 
By the symmetry between $u$ and $v$, it is enough to prove that there is a matching in $\ecol$ which covers $u$ and misses $v$. Suppose to the contrary that every matching in $\ecol$ covering $u$ covers $v$ too. By i), every $M\in\ecol$ covers $v$, so $\chi'(H-e)=d_{H}(v)-1\leq \Delta(H)-1$. Thus $\chi'(H)\leq \Delta(H)$, which contradicts the assumption on $H$.
\end{dem}
 
For every graph $H$ and each $F\se\eds{H}$, we write $\psbg{H}{F}$ for the graph $(\vts{H},F)$ (ambiguity with the notation for subgraphs induced by sets of vertices should not occur). We will often use the following basic recoloring-argument, which corresponds to switching the colors on a bi-colored component of a graph.

\begin{prop}\label{prop-edge-recol}
Let $H$ be a graph, $\ecol$ be an edge-coloring of $H$ and $A,B$ be distinct elements of $\ecol$.
If $K$ is a component of $\psbg{H}{A\Delta B}$, then $ (\ecol\sm\set{A,B})\cup\set{A\Delta \eds{K},B\Delta \eds{K} } $ is an edge-coloring of $H$ which uses $|\ecol|$ colors. 
\end{prop}
We give a few more notations for the proof of \myrefwp{thm-ben-crit}. 

Let $G$ and $H$ be graphs. 
The graph $G\cap H$ is defined by $\vts{G\cap H}=\vts{G}\cap\vts{H}$ and $\eds{G\cap H}=\eds{G}\cap\eds{H}$. 

Let $R$ be an odd ring of $H$. A matching $M$ of $H$ is an \emph{$R$-matching} if  $|\eds{R}\cap M|=\frac{|\vts{R}|-1}{2}$ and $M$ misses a (necessarily unique) vertex of $C$.
A \emph{chord} of a circuit $C$ is an edge between two non-adjacent vertices of $C$. The \emph{ends} of a path $P$ are its (at most 2) vertices of degree less than 2. The \emph{end-edges} of $P$ are the edges of $P$ (if any) which are incident to its ends.

We shall detect odd-$C_5^{+}$ subgraphs of $H$ using the following basic remark:  \emph{the odd-$C_5^{+}$ subgraphs of $H$ are the simple graphs formed by an odd circuit $C$ and an odd path $P$ of $H$ such that $\vts{C}\cap\vts{P}$ is reduced to the ends of $P$.} In particular, an odd circuit with a chord forms an odd-$C_5^{+}$.

\begin{dem}[of lemma \ref{thm-ben-crit}]
Let $H$ be a graph (with possibly multiple edges) such that $\chi'(H)>\Delta(H)$ and let $e$ be a critical edge of $G$ which is not an edge of an odd-$C_5^{+}$ of $H$. 

Let $u$ and $v$ be the ends of $e$ and $\ecol$ be an optimal edge-coloring of $H-e$. By \myrefwp{lem-precol},
there exist matchings $A,B\in\ecol$ such that $A$ covers $u$ and misses $v$, whereas $B$ covers $v$ and misses $u$.

Consider the component $P$ of $u$ in the graph $\psbg{H}{A\Delta B}$. It is either a path or a circuit, but $B$ misses $u$ so it must be a path. We have
$ v\in\vts{P} $:
otherwise, 
$(\ecol\sm\set{A,B})
\cup\set{A\Delta \eds{P},B\Delta\eds{P}}$ 
would be an optimal edge-coloring of $H-e$ (by proposition \ref{prop-edge-recol}) in which $A\Delta\eds{P}$ misses both $u$ and $v$. It could therefore be extended to an edge-coloring of $H$, contradicting $\chi'(H-e)<\chi'(H)$.

So $P$ is a $uv$-path and the circuit $L$ of $H$ obtained by adding $e$ to $P$ is odd. If $L$ had a chord $f$, then $L$ and $f$ would form an odd-$C_5^{+}$ of $H$ containing $e$: a contradiction.

Thus, $L$ is an induced circuit and $\vts{L}$ induces an odd ring $R$ of $H$. Let $r=\frac{|\vts{R}|-1}{2}$. We claim that:
\begin{equation}\label{crit-claim-1}
\text{\emph{$R$ contains exactly $r$ edges of each matching $M$ of $\ecol$.}}
\end{equation}
Let us immediately show that this claim implies the theorem: excepting $e$, every edge of $R$  belongs to a matching of $\ecol$ and these matchings are pairwise-disjoint. Therefore, the number of edges of $R$ is $|\ecol| r+1=r\cdot\chi'(H-e)+1$ and the conclusion follows.

We now prove \eqref{crit-claim-1}.
Let $M\in\ecol$. 
If $M\in\set{A,B}$, then $M$ has $r$ edges in $R$ because the edges of $P$ alternate  between $A$ and $B$.
So let us henceforth assume that $M\notin\set{A,B}$. 
Using the symmetry between $u$ and $v$, we may suppose without loss of generality that $M$ covers $u$.
Let $K$ be the component of $u$ in $\psbg{H}{M\Delta B}$. The graph $K$ is a path since $B$ misses $u$. We have:
\begin{equation}\label{crit-claim-2}
\text{\emph{$K\cap R$ is an even path.}}
\end{equation}
Suppose to the contrary that $K\cap R$ has more than one component. Then, there exists a (non-zero length) path $Q$ of $K$ which is edge-disjoint from $R$ and whose ends belong to $\vts{R}$. Since $B$ is an $R$-matching, both end-edges of $Q$ must belong to $M$ and $Q$ is odd.
Hence, the graph $L$ and $Q$ together form an odd-$C_5^{+}$ of $H$ containing $e$: a contradiction with our assumption. So $K\cap R$ is connected and it joins $u$ to some vertex $w$ of $R$. 

Suppose $|\eds{K\cap R}|>0$. Since $B$ is an $R$-matching missing $u$, the vertex $w$ is covered by $B$ in the graph $K\cap R$. Therefore, the path $K\cap R$ has exactly one end-edge in $B$ and it must be even as stated above.

Let $M'=M\Delta \eds{K}$. By \eqref{crit-claim-2} and since the edges of $K\cap R$ alternate between $M$ and $B$, we obtain:
 \[ |M\cap \eds{R}|=|M'\cap \eds{R}|, \]

We now show that $|M'\cap \eds{R}|=r$. This will end the proof of \eqref{crit-claim-1}.

Let $B'=B\Delta \eds{K}$  and $\ecol'=(\ecol\sm\set{M,B})\cup\set{M',B'}$. By proposition \ref{prop-edge-recol}, $\ecol'$ is an optimal edge-coloring of $H-e$. Notice that $M'$ misses $u$. Therefore, $M'$ must cover $v$:  otherwise, $\ecol'$ could be extended to an  edge-coloring of $H$ by adding $e$ to $M'$ and this would contradict $\chi'(H-e)<\chi'(H)$.

Let $K'$ be the component of $v$ in $\psbg{H}{M'\Delta A}$ and $T=K'\cap R$. Since $A$ is an $R$-matching of $H$ which misses $v$,  
we can repeat the argument of the proof of \eqref{crit-claim-2} to show that $T$ is a path.
Now, we have:
\[ u\in \vts{T}.  \]
Indeed, suppose that $u\notin\vts{T}$ and let $M''=M\Delta \eds{K'}$, $A'=A\Delta \eds{K'}$ and $\ecol''=(\ecol'\sm\set{M',A})\cup\set{M'',A'}$. 
As above, proposition \ref{prop-edge-recol} implies that  $\ecol''$ is an optimal edge-coloring of $H-e$. However, $M''$ misses both ends of $e$ so $\ecol''$ can be extended to an edge-coloring of $H$ by adding $e$ to $M''$: a contradiction.

Since $M'$ does not contain an edge parallel to $e$ (it misses $u$) and since $u\in\vts{T}$, the only way for $T$ to be a path is that it coincides with $P$ in the underlying simple graph of $H$. But $T$ alternates between $M'$ and $A$, hence $|M'\cap \eds{R}|=|M'\cap \eds{R}|=r$.

\end{dem}

\section{T-perfect claw-free graphs}
Our purpose is to prove theorem \ref{thm-ben-cft-irp}. The first part gives an outline of the proof and hopefully clarifies that we have to take a new approach compared to the unweighted case. The proofs of the two main lemmas are postponed to sections 4.2 and 4.3.

\subsection{How the proof works ?}
The unweighted case of theorem \ref{thm-ben-cft-irp} was obtained by Bruhn and Stein in \cite{Bruhn2012} and appears as a preliminary result of theorem \ref{thm-brst-ir}. It is not difficult to see that it means that t-perfect claw-free graphs are 3-colorable.
\begin{thm}[Bruhn, Stein]\label{thm-bs-3col}
Every  t-perfect claw-free graph is 3-colorable.
\end{thm}
Let us briefly recall the approach of the proof: the result is first proved for line graphs (using edge-colorings). Then, considering a t-perfect claw-free graph they show: \emph{if $G$ is 3-connected then it is either a line graph or one of a few exceptional graphs which can be easily 3-colored. }
Otherwise, they use a 2-cut (a subset of at most two vertices of $G$ whose deletion disconnects $G$) to decompose the graph into smaller pieces and apply induction (which is not straightforward).

Unfortunately, it is not straightforward how two weighted colorings with a small number of colors can be combined along a 2-cut such that the number of colors remains small. Kilakos and Marcotte \cite{Kilakos1997} gave general sufficient conditions under which this operation can be performed. However, it is not clear how to apply it directly to t-perfect claw-free graphs. 
Therefore, we follow a different approach.

We proceed by an induction where the line graphs form the base case. In the presence of certain subgraphs, we reduce the weight function. If no such subgraph appears, then we show that the graph considered is a line graph and apply theorem \ref{thm-ben-line-irp}.

\begin{figure}
\centering
\includegraphics{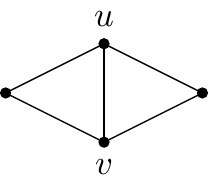}
\caption{A diamond with central vertices $u,v$}\label{fig-diamond}
\end{figure}

A \emph{diamond} of a graph $G$ is an induced subgraph $D$ of $G$ which is isomorphic to the complete graph $K_4$ minus an edge (see figure \ref{fig-diamond}). A vertex $u$ of $D$ is \emph{central} if $d_D(u)=3$. 

A central vertex $u$ of $D$ is \emph{small} if $d_G(u)=3$.
We say that $D$ is \emph{small} if it has a small central vertex and that it is \emph{large} otherwise. Notice that if $D$ is large, then both of its central vertices have degree at least 4 in $G$.

\begin{lem}\label{lem-reduc}
Let $G$ be a t-perfect claw-free graph and $c\in\zset_+^{\vts{G}}$.
Suppose that $G$ has a small diamond $D$ with a small central vertex $v$ such that $c_v\geq 1$.
Put $c'=c-\chi^{v}$. 

If $\cn{c'}{G}=\ceil{\fcn{c'}{G}}$, then $\cn{c}{G}=\ceil{\fcn{c}{G}}$.
\end{lem}

So an induction on the size of the weight function can be performed as long as there is a small diamond in the graph. 
Now, the following result shows that the remaining case falls in the scope of theorem \ref{thm-ben-line-irp}.

\begin{lem}\label{lem-small-diamond}
Let $G$ be a t-perfect claw-free simple graph. If every diamond of $G$ is large, then $G$ is a line graph.
\end{lem}
The starting point of the approach of Bruhn and Stein for theorem \ref{thm-brst-ir} is the use of Harary's characterization of line graphs of simple graphs (in terms of triangles). It plays a similar role in the proof of lemma \ref{lem-small-diamond}.
The other key-ingredient is the characterization of t-perfection among squares of circuits. These two results are stated in part 4.3.

We now prove that theorem \ref{thm-ben-cft-irp} follows from lemmas \ref{lem-reduc}, \ref{lem-small-diamond} and theorem \ref{thm-ben-line-irp}.

\begin{dem}[of theorem \ref{thm-ben-cft-irp}]
Seeking a contradiction, consider a counter-example $(G,c)$ which is minimum with respect to $|\vts{G}|+\sum_{v\in\vts{G}}c_v$. Clearly, we can assume that $G$ is simple and that no coordinate of $c$ is equal to zero. Thus, $G$ cannot have a small diamond because of lemma \ref{lem-reduc}. Therefore, lemma \ref{lem-small-diamond} shows that $G$ is a line graph and the conclusion follows from theorem \ref{thm-ben-line-irp}.
\end{dem}

\subsection{Proof of lemma \ref{lem-reduc}}
\begin{figure}
	\centering
	\subfloat[$N=3$]
		{\includegraphics{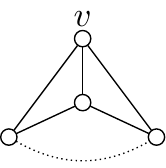}\label{fig-n3}}\qquad
	\subfloat[$N=4$]
		{\includegraphics{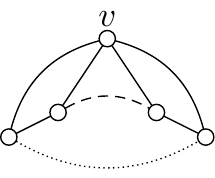}\label{fig-n4}}\qquad
	\subfloat[$N=5$]
		{\includegraphics{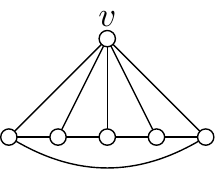}\label{fig-n5}}
	\caption{The different possibilities for $H$: 
	{\small dotted and dashed lines denote pairwise-disjoint paths. Each ordinary line denotes an edge,  dotted lines correspond to odd paths and the dashed one to an even path of non-zero length. There is no other edge.}}\label{fig-prop-neighbor}
\end{figure}

Let $G$ be a graph and $v\in\vts{G}$. We write $N_G(v)$ for the set of neighbors of $v$. In particular, if $G$ is simple then $d_G(v)=|N_G(v)|$ and $\Delta(G)$ is the largest number of neighbors of a vertex of $G$.
We will need a basic property of odd holes in t-perfect claw-free graphs:

\begin{prop}\label{prop-neighbor}
Let $G$ be a t-perfect claw-free graph. If $v\in\vts{G}$ and $C$ is an induced odd-circuit of $G$, then $v$ has at most 2 neighbors in $C$. 
\end{prop}
\begin{dem}
Let $N=|N_G(v)\cap\vts{C}|$ and $H=G\cro{\vts{C}\cup\set{v}}$.  To the contrary, suppose that $N\geq 3$.
Since $G$ is claw-free, $N\leq 5$ and $H$ is either of the form shown in figures \ref{fig-n3} and \ref{fig-n4} (if $N\leq 4$) or it is isomorphic to the graph of figure \ref{fig-n5} (if $N=5$). 

However, none of these graphs are t-perfect. Indeed, 
let $k=\frac{|\vts{C}|-1}{2}$ and $y$ be the vector of $\qset^{\vts{H}}$ whose coordinates are equal to $\frac{k}{|\vts{C}|}$ on every vertex of $C$ and to $\frac{1}{|\vts{C}|}$ elsewhere. Clearly, $y$ satisfies the non-negativity, clique and odd-hole inequalities of $H$ but it does not belong to $\ssp{H}$: the maximum over $\ssp{H}$ of the linear function $\sum_{u\in\vts{H}}x_u$ is $k$ (it is the maximum cardinality of a stable set of $H$) but $\sum_{u\in\vts{H}}y_u>k$.
Therefore, $H$ is not t-perfect. By proposition \ref{prop-tperf-indu}, this contradicts the t-perfection of $G$.
\end{dem}

We now prove lemma \ref{lem-reduc}. We will use again color exchanges on bi-colored components but for colorings of the vertices (the statement of proposition \ref{prop-edge-recol} can be easily translated using line graphs).

\begin{dem}[of lemma \ref{lem-reduc}]
We start with an argument similar to the proof of lemma \ref{lem-precol}.
Let $D=G\cro{\set{x,v,w,y}}$, where $x$ and $y$ are the two vertices of degree 2 in $D$. Hence, the neighbors of $v$ in $G$ are $x$, $y$ and $w$.
Let $\mathcal{F}$ be a coloring of $(G,c')$. Without loss of generality, we can assume that every $u\in\vts{G}$ belongs to exactly $c'_u$ members of $\mathcal{F}$.

If there exists an $S\in\mathcal{F}$ such that $S\cap(N_G(v)\cup\set{v})=\vn$, then $(\mathcal{F}\sm\set{S})\cup\set{S\cup\set{v}}$ is a coloring of $(G,c)$ with $\cn{c'}{G}$ colors. Clearly, $\fcn{c'}{G}\leq\fcn{c}{G}$ and the result of the lemma follows. 

Hence we will assume that every member of $\mathcal{F}$ meets $ N_G(v)\cup\set{v} $.
For every $u\in\vts{G}$, let $\mathcal{F}_u$ denote the set of members of $\mathcal{F}$ containing $u$.
First, suppose that $\mathcal{F}_x\se\mathcal{F}_y$. Then, the number of members of $\mathcal{F}$ intersecting $ N_G(v)\cup\set{v} $ is:
\[ \chi(G,c')=|\mathcal{F}_x\cup\mathcal{F}_y\cup\mathcal{F}_w\cup\mathcal{F}_v|
=|\mathcal{F}_y\cup\mathcal{F}_w\cup\mathcal{F}_v|
=|\mathcal{F}_y|+|\mathcal{F}_w|+|\mathcal{F}_v|\leq\omega(G,c)-1\leq\ceil{\fcn{c}{G}}-1, \]
as $\set{v,w,y}$ is a clique (the last inequality follows from proposition \ref{prop-fcn-hparf}). 
So adding $\set{v}$ to $\mathcal{F}$ gives a coloring of $(G,c)$ with $\ceil{\fcn{c}{G}}$ colors and we are done.

Therefore, we may assume that $\mathcal{F}_x\nsubseteq\mathcal{F}_y$ and by symmetry that $\mathcal{F}_y\nsubseteq\mathcal{F}_x$.
Let $S\in \mathcal{F}_x\sm\mathcal{F}_y$ and $T\in \mathcal{F}_y\sm\mathcal{F}_x$. 
Consider $H=G\cro{S\Delta T}$. 
Let $K$ be the component of $x$ in $H$. 
We claim that:
\begin{center}
$y\notin \vts{K}.$
\end{center}
Suppose to the contrary that $y\in\vts{K}$ and
let $P$ be a shortest (thus induced) path of $K$ joining $x$ and $y$.
The vertices of $P$ alternate between $S$ and $T$, hence $P$ has odd length. As $G$ does not have a clique of size 4 (it is t-perfect), the length of $P$ is at least 3. 
Let $L=G\cro{\vts{P}\cup\set{v}}$. 
By assumption, $d_G(v)=3$ so $L$ is an odd hole of $G$. 
Now, $w$ does not belong to $L$ (as it is adjacent to $x$ and $y$) and it has at least 3 neighbors in $L$. By proposition \ref{prop-neighbor}, this contradicts the t-perfection of $G$.

The lemma now easily follows: $(\mathcal{F}\sm\set{S,T})\cup\set{S\Delta \eds{K}, T\Delta \eds{K}}$ is a coloring of $(G,c)$ with $\chi(G,c')$ colors and $\cn{c'}{G}=\ceil{\fcn{c'}{G}}\leq\ceil{\fcn{c}{G}}$.

\end{dem}

\subsection{Proof of lemma \ref{lem-small-diamond}}
We now consider the family of squares of circuits. 
Let $n\geq 3$ be an integer. The graph $C_n^{2}$ is defined as follows: $\vts{C_n^{2}}=\set{1,2,\ldots,n}$ and $\eds{C_n^{2}}$ is the set of pairs $\set{i,j}$ of integers between 1 and $n$ such that $|i-j|\leq 2$ (mod $n$). The following theorem is a consequence of results of Dahl \cite{Dahl1999} and was obtained in a rather different way by Bruhn and Stein \cite{Bruhn2012} using \emph{t-minors}. 

\begin{figure}[h!]
\centering
\includegraphics{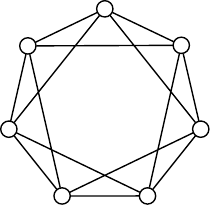}
\caption{The graph $C_7^{2}$}\label{fig-c72}
\end{figure}

There are two t-minor operations: the deletion of a vertex, and the contraction of every edge incident to a vertex whose neighborhood is a stable set. Gerards and Shepherd \cite{Gerards1998} proved that \emph{these transformations keep t-perfection}.
In \cite{Bruhn2012}, it is shown that \emph{$C_7^{2}$ \emph{(see figure \ref{fig-c72})} and $C_{10}^{2}$ are minimally non-t-perfect} with respect to t-minor operations, whereas $K_4$ (which is not t-perfect) is a t-minor of $C_n^{2}$ for every $n\geq 5$ which is not 6, 7 or 10.
We need only a subset of these results:

\begin{thm}[Bruhn, Stein \cite{Bruhn2012}]\label{lem-brst-squares}
Let $n\geq 3$ be an integer. The graph $C_n^{2}$ is  t-perfect if and only if $n\in\set{3,6}$.
\end{thm}

We will use that t-perfect claw-free simple graphs have small degree. This already appeared in \cite{Bruhn2012} but we reproduce the proof here for the commodity of the reader.

\begin{prop}\label{prop-degree-tpc}
Every t-perfect claw-free simple graph has maximum degree at most 4.
\end{prop}
\begin{dem}
Let $G$ be a t-perfect claw-free simple graph and let $v\in\vts{G}$. Since $G$ is simple, we have $d_G(v)=|N_G(v)|$.

First, notice that $d_G(v)\leq 5$: otherwise, by Ramsey's theorem, $N_G(v)$ would contain either a triangle or a stable set of cardinality 3. 
So $G$ would contain a clique of size 4 or a claw: a contradiction.

Now, if $d_G(v)=5$ then $G\cro{N_G(v)}$ is a graph with 5 vertices having no stable set of size 3 and no triangle. Hence it is an odd circuit of length 5 and this contradicts \myrefwp{prop-neighbor} (in part 4.2).
\end{dem}

The last ingredient needed to prove lemma \ref{lem-small-diamond} is the following proposition on claw-free simple graphs.
Recall that $\omega(G)$ denotes the maximum cardinality of a clique of a graph $G$.
\begin{prop}\label{prop-large-diamonds}
Let $G$ be a connected claw-free simple graph such that $\Delta(G)\leq 4$ and $\omega(G)\leq 3$. If every diamond of $G$ is large, then at least one of the following statements holds:
\begin{itemize}
\item[i)] $G$ is a line graph,
\item[ii)] there exists an integer $k\geq 7$ such that $G$ is isomorphic to $C_k^{2}$.
\end{itemize}
\end{prop}
The proof of this proposition is postponed to the end of this section. We first show that these results imply lemma \ref{lem-small-diamond}:

\begin{dem}[of lemma \ref{lem-small-diamond}]
Let $G$ be a t-perfect claw-free simple graph such that every diamond of $G$ is large. Clearly, we need only to prove that each component of $G$ is a line graph. 

Let $H$ be a component of $G$. Since $H$ is an induced subgraph of $G$, we have that $H$ is claw-free and by proposition \ref{prop-tperf-indu}, it is also t-perfect.
Hence, using proposition \ref{prop-degree-tpc} we obtain that $\Delta(H)\leq 4$. 
Furthermore, $\omega(H)\leq 3$ and every diamond of $H$ is large.

Suppose to the contrary that $H$ is not a line graph. By proposition \ref{prop-large-diamonds}, there exists an integer $k\geq 7$ such that $H$ is isomorphic to $C_k^{2}$. Now, theorem \ref{lem-brst-squares} shows that $H$ is not t-perfect. By proposition \ref{prop-tperf-indu}, this contradicts the t-perfection of $G$.

\end{dem}

We end this section with the proof of proposition \ref{prop-large-diamonds}. We need a characterization of line graphs by Harary in terms of diamonds.
A triangle $T$ of a graph $G$ is \emph{odd} if $G$ contains a vertex $v\notin T$ which has an odd number of neighbors in $T$. A diamond of $G$ is \emph{odd} if both of its triangles are odd triangles of $G$. 
The implication ii)$\Rightarrow$i) of the following result is the key to obtain line graphs:

\begin{thm}[Harary \cite{Harary1972}, pg. 74: theorem 8.4]\label{thm-harary}
Let $G$ be a claw-free simple graph. The following statements are equivalent:
\begin{itemize}
\item [i)] $G$ is the line graph of a simple graph,
\item [ii)] $G$ does not have an odd diamond.
\end{itemize}
\end{thm}
Actually, we do not use the simplicity of the graph $H$ whose line graph is $G$ since theorem \ref{thm-ben-line-irp} holds for line graphs of non-necessarily simple graphs.

\begin{dem}[of proposition \ref{prop-large-diamonds}]
Let $G$ be a connected claw-free simple graph with $\Delta(G)\leq 4$, $\omega(G)\leq 3$ and such that every diamond of $G$ is large. 
Furthermore, let us assume that $G$ is not a line graph. 
We have to prove:
 
\begin{center}
\emph{there exists an integer $k\geq 7$ such that $G$ is isomorphic to the graph $C_k^{2}$.}
\end{center}

By theorem \ref{thm-harary}, $G$ has an odd diamond $D$. 
Put $D=G\cro{v_1,v_2,v_3,v_4}$, where $v_2$ and $v_3$ are the central vertices of $D$. 
Since $D$ is large, we have $d_G(v_2)=4$. Hence, $v_2$ has a neighbor $v_5\in\vts{G}\sm\set{v_1,v_3,v_4}$. As $G$ is claw-free, $v_5$ is adjacent to at least one of $v_1$ and $v_4$. Using the symmetry between $v_1$ and $v_4$, we can assume without loss of generality that $v_4v_5\in\eds{G}$.

Again, $D$ is large so $ d_G(v_3)=4 $. Since $v_5$ cannot be a neighbor of $v_3$ (otherwise $\set{v_2,v_3,v_4,v_5}$ would be a clique of size 4), there exists $v_6\in\vts{G}\sm\set{v_1,\ldots,v_5}$ such that $v_3v_6\in\eds{G}$. But $G$ is claw-free so at least one of $v_1v_6$ and $v_4v_6$ is an edge of $G$. However:
\begin{equation}\label{claim-diamonds-1}
v_4v_6\notin\eds{G}.
\end{equation}
Otherwise, $v_1$, $v_5$ and $v_6$ would be the only vertices of $G$ having a neighbor in the triangle $v_2v_3v_4$ (because $\Delta(G)=4$).
But each of these vertices have exactly two neighbors on $v_2v_3v_4$. Thus $v_2v_3v_4$ would not be an odd triangle of $G$. 
This contradicts that $D$ is odd.

Therefore, $v_1v_6\in\eds{G}$. Since $v_1$ is a central vertex of the diamond induced by $\set{v_1,v_2,v_3,v_6}$ (and every diamond of $G$ is large), we must have $d_G(v_1)=4$. Furthermore:
\begin{equation}\label{claim-diamonds-2}
v_1v_5\notin\eds{G}.
\end{equation}
Else, the same argument used to prove \eqref{claim-diamonds-1} shows that $v_1v_2v_3$ would not be an odd triangle, contradicting that $D$ is odd.

So $v_1$ must have a neighbor $v_7\notin\set{v_2,\ldots,v_6}$ and as $G$ is claw-free, $v_7$ is adjacent to at least one of $v_2$ and $v_6$. But $v_2$ already has 4 neighbors among $\set{v_1,v_3,v_4,v_5,v_6}$, thus $v_2v_7\notin\eds{G}$ and $v_6v_7\in\eds{G}$. 

Let $H=G\cro{\set{v_1,\ldots,v_7}}$ (see figure \ref{fig-proof-dia1}). 
\begin{figure}[h]
\centering
\includegraphics[scale=1.5]{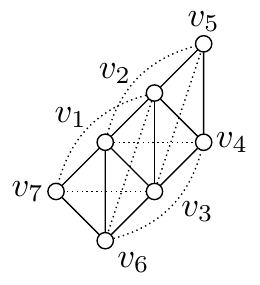}
\caption{{\small The construction of $v_1,\ldots,v_7$. 
Dotted lines indicate pairs of non-adjacent vertices.}}\label{fig-proof-dia1}
\end{figure}

\begin{case}
\emph{$N_G(v_4)$ and $N_G(v_6)$ are contained in $\vts{H}$}. 

Since $v_4$ and $v_6$ are both central vertices of diamonds of $G$, we have $d_G(v_6)=d_G(v_4)=4$. Recall that $v_4v_6\notin\eds{G}$. Thus,  $v_5v_6\in\eds{G}$ and $v_4v_7\in\eds{G}$. In particular, the vertices $v_3$, $v_5$ and $v_7$ are neighbors of $v_4$. But $G$ has no clique of size 4 so $v_3v_7\notin\eds{G}$ and $v_3v_5\notin\eds{G}$. Since $G$ is claw-free, this implies that  $v_5v_7\in\eds{G}$ (see figure \ref{fig-proof-dia2}).

Now, the map $1\rightarrow v_7$, $2\rightarrow v_6$, $3\rightarrow v_1$, $4\rightarrow v_3$, $5\rightarrow v_2$, $6\rightarrow v_4$, $7\rightarrow v_5$ defines an isomorphism from $C_7^{2}$ to a subgraph of $H$. 
Since $C_7^{2}$ is 4-regular and $\Delta(G)=4$, the graph $H$ is in fact isomorphic to $C_7^{2}$ and is a component of $G$. Since $G$ is connected, $G=H$ and the conclusion follows.
\end{case}

\begin{figure}[h]
\centering
\includegraphics[scale=1.5]{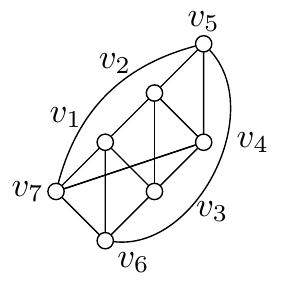}
\caption{{\small \textbf{case 1}: building a $C_7^{2}$ from $v_1,\ldots,v_7$ when $N_G(v_4)\cup N_G(v_6)\se\vts{H}$. }}\label{fig-proof-dia2}
\end{figure}

\begin{figure}[h!]
\centering
\includegraphics[scale=1.5]{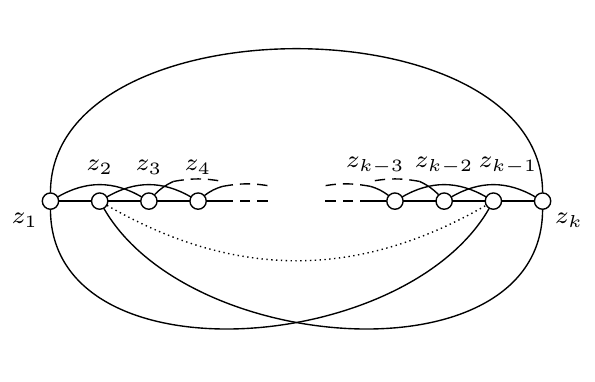}
\caption{{\small \textbf{case 2}: building an induced $C_k^{2}$ (with $k\geq 8$) from $Q$.
The dotted line indicates that $z_2z_{k-1}\notin\eds{G}$}}\label{fig-proof-dia3}
\end{figure}
\begin{case}
\emph{At least one of $v_4$ and $v_6$ has a neighbor outside of $\vts{H}$.}

Using the symmetry between $v_4$ and $v_6$, we can assume without loss of generality that $v_6$ has a neighbor $v_8\notin \vts{H}$. Now, the vertices $v_3$, $v_7$ and $v_8$ are three neighbors of $v_6$. Since $G$ has no clique of size 4, we have $v_3v_7\notin\eds{G}$. Furthermore, $v_3$ already has 4 neighbors among the vertices of $H$ so $v_3v_8\notin\eds{G}$ and since $G$ is claw-free we have $v_7v_8\in\eds{G}$. Define $w_1=v_8$, $w_2=v_7$, $w_3=v_6$, $w_4=v_1$ $w_5=v_3$, $w_6=v_2$, $w_7=v_4$, $w_8=v_5$. Notice that $(w_1,\ldots,w_8)$ is a path such that $w_iw_{i+2}\in\eds{G}$ for every $1\leq i\leq 6$. 

Let $Q=(z_1,\ldots,z_k)$ be a path of $G$ such that $z_iz_{i+2}\in\eds{G}$ for every $1\leq i\leq k-2$, and choose
$k$ maximum. The path $(w_1,\ldots,w_8)$ shows that $k\geq 8$. 
Let $L=G\cro{z_1,\ldots,z_k}$. We claim that:
\begin{equation}\label{claim-diamonds-3}
\text{\emph{$L$ is isomorphic to $C_k^{2}$.}}
\end{equation}
Since $\Delta(G)=4$ and $C_k^{2}$ is 4-regular, this implies that $L$ is a component of $G$. As $G$ is connected, we have $G=L$ and this ends the proof of proposition \ref{prop-large-diamonds}. 
We now prove \eqref{claim-diamonds-3}.

Since $z_2$ is a central vertex of the diamond induced by $\set{z_1,z_2,z_3,z_4}$, we have $d_G(z_2)=4$. 
Furthermore, $N_G(z_2)\se\vts{L}$: we could otherwise use a vertex $z_{0}\in N_G(z_2)\setminus\vts{L}$ to extend $Q$ ($z_0$ must be adjacent to $z_1$ because $G$ is claw-free) and contradict the maximality of $k$. 
For the same reason, $N_G(z_{k-1})\se\vts{L}$.

Now, notice that for every $i\in\set{3,\ldots,k-2}$: $d_L(z_i)=4$ .
Hence, at least one of $z_k$ and $z_{k-1}$ is a neighbor of $z_2$.
However, $z_2z_{k-1}$ cannot be an edge of $G$: 
the vertices $z_1$, $z_4$ and $z_{k-1}$ would otherwise be three pairwise non-adjacent neighbors of $z_2$ and $G$ would contain an induced claw: a contradiction.

Therefore, $z_2z_{k}\in\eds{G}$. Similarly, $d_G(z_{k-1})=4$ and we have $z_1z_{k-1}\in\eds{G}$. 
Since $G$ is claw-free, this implies that $z_1z_k\in\eds{G}$ (see figure \ref{fig-proof-dia3}). Finally, the map $i\rightarrow z_i$ defines an isomorphism from $C_k^{2}$ to a subgraph of $L$. Since $C_k^{2}$ is 4-regular, $L$ is in fact isomorphic to $C_k^{2}$.

\end{case}

\end{dem}

\section{Minmax formulae and algorithmic remarks}
In this section, we first obtain an explicit formula for the weighted chromatic number of t-perfect claw-free graphs and h-perfect line-graphs. We also give a corresponding formula for the chromatic index of an \cfpfree graph without referring to its line graph. Finally, we discuss the algorithmic aspects of our results.

Using proposition \ref{prop-fcn-hparf} and theorems \ref{thm-ben-line-irp} and \ref{thm-ben-cft-irp}, we  obtain:

\begin{cor}
Let $G$ be a graph and $c\in\wght{G}$. If $G$ is a t-perfect claw-free graph or an h-perfect line-graph, then:
\[  \cn{c}{G}=\max(\omega(G,c),\ceil{\Gamma(G,c)}). \]
\end{cor}

It is easy to check that $\Gamma(G,c)\leq 3$ when $c$ is the all-1 vector, thus we obtain the 3-coloring result of Bruhn and Stein (theorem \ref{thm-bs-3col}) in another way.

By a theorem of Edmonds and Pulleyblank \cite{Pulleyblank1974}, \emph{the non-trivial facets of the matching polytope of a graph are given by the inequalities defined by stars and 2-connected factor-critical induced subgraphs}. 
A consequence of this theorem (through a result of Lov\'asz \cite{Lovasz1972a} which relates factor-criticality and ear-decompositions) is that \emph{for every \cfpfree graph $H$: $\chi_f'(H)=\max(\Delta(H),\Gamma'(H))$}. We do not present the derivation of this formula because it essentially amounts to repeat the proof of the characterization of h-perfect line-graphs in terms of odd-$C_5^{+}$ subgraphs, which can be found in \cite{Cao1998} and \cite{Bruhn2012}. Using theorem \ref{thm-ecol-form}, we obtain:

\begin{cor}\label{cor-cfpfree}
Every \cfpfree graph $H$ satisfies $\chi'(H)=\max(\Delta(H),\ceil{\Gamma'(H)})$.
\end{cor}

\begin{figure}[h]
\centering
\includegraphics{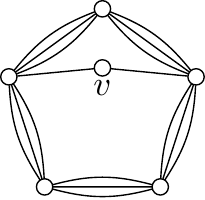}
\caption{The graph $H_3$}\label{fig-tightness}
\end{figure}

The difference $\Gamma'-\Delta$ can be arbitrarily large for \cfpfree graphs. Indeed, let $m$ be a positive integer and let the graph $H_m$ be obtained as follows (see also figure \ref{fig-tightness}): start with a circuit $C$ of length 5, replace every edge of $C$ with $m$ parallel edges and add a new vertex $v\notin\vts{C}$ adjacent to exactly two non-adjacent vertices of $C$. Clearly, $H_m$ is \cfpfree and $\Delta(H_m)=2m+1$ whereas $\Gamma'(H_m)=\frac{5m}{2}$.

As a first algorithmic remark, notice that it is straightforward to turn the proofs of proposition \ref{lem-precol} and lemma \ref{thm-ben-crit} into a polynomial-time algorithm which computes an optimal edge-coloring of an \cfpfree graph.

The \emph{maximum-weight stable set problem} can be formulated as follows: considering a graph $G$ and a weight $c\in\wght{G}$, find a stable set $S$ of $G$ such that $\sum_{v\in S}c_v$ is maximum. Gr\"{o}tschel, Lov\'asz and Schrijver \cite{Grotschel1986} proved that \emph{this problem is solvable in polynomial-time for h-perfect graphs}. Since we only consider claw-free graphs, we can use a specific algorithm (whose construction is more elementary):

\begin{thm}[Minty \cite{Minty1980}, Sbihi \cite{Sbihi1980}, see also Nakamura and Tamura \cite{Nakamura2001}]
The maximum-weight stable set problem can be solved in polynomial-time in the class of claw-free graphs.
\end{thm} 
By results of Gr\"{o}tschel, Lov\'asz and Schrijver \cite{Groetschel1988}, this implies that 
\emph{there exists a polynomial-time algorithm to find the weighted fractional chromatic number of $(G,c)$ for every claw-free graph $G$ and every $c\in\wght{G}$.}
Adding a rounding-step to this algorithm and using theorems \ref{thm-ben-line-irp} and \ref{thm-ben-cft-irp}, we obtain:
\begin{cor}
There exists a polynomial-time algorithm which computes $\cn{c}{G}$ for every graph $G$ and every weight $c\in\wght{G}$ such that $G$ is either an h-perfect line-graph or a  t-perfect  claw-free graph.
\end{cor}

We do not know of a \emph{combinatorial} polynomial-time algorithm which computes the fractional chromatic number of a t-perfect graph. 
Furthermore, contrarily to the case of perfect graphs, we cannot directly derive the existence of an efficient algorithm to find an \emph{optimal coloring} of the graph from our results.

\bibliographystyle{alpha}
\bibliography{idpbib}

\end{document}